\documentclass[12pt]{report}
\usepackage{amssymb,amsmath,makeidx,verbatim}
\newcommand{\C}{\mathbb{C}}

\newcommand{\be}{\begin{enumerate}}
	\newcommand{\ee}{\end{enumerate}}
\newcommand{\bq}{\begin{eqnarray*}}
	\newcommand{\eq}{\end{eqnarray*}}

\begin{document}
	%\pagenumbering{roman}
	\newcommand{\disp}{\displaystyle}
	\thispagestyle{empty}
	\begin{center}
		\textsc{A note on the $L^{2}-$harmonic analysis of the Joint-Eigenspace Fourier transform\\}
		\ \\
		\textsc{Olufemi O. Oyadare}\\
		\ \\
		Department of Mathematics,\\
		Obafemi Awolowo University,\\
		Ile-Ife, $220005,$ NIGERIA.\\
		\text{E-mail: \textit{femi\_oya@yahoo.com}}\\
	\end{center}
	\begin{quote}
		{\bf Abstract.} {\it We consider the irreducibility of the regular representation of a noncompact semisimpe Lie group $G$ on the Hilbert space of the image of the Joint-Eigenspace Fourier transform on its corresponding symmetric space $G/K.$ The $L^{2}-$decomposition of the Joint-Eigenspace Fourier transform leads to the complete characterization of the said irreducibility in terms of the simplicity of a pair of members of $\mathfrak{a}^{*}_{\mathbb{C}}.$}
	\end{quote}
	
	\ \\
	\ \\
	\ \\
	\ \\
	\ \\
	\ \\
	\ \\
	\ \\
	\ \\
	\ \\ 
	\ \\
	\ \\
	\ \\
	\ \\
	\ \\
	\ \\
	\ \\
	$\overline{2020\; \textmd{Mathematics}}$ Subject Classification: $53C35, \;\; 43A90, \;\; 42A38$
	
	\ \\
	{\bf \S1. Introduction.}
	
	\indent The present paper is a continuation of results on the complete Fourier transform on a noncompact symmetric space $X=G/K$ (for a maximal compact subgroup $K$ of a noncompact semisimple Lie group $G$ with finite centre $[3.]$) introduced as the Joint-Eigenspace Fourier transform on $X$ in the paper of Oyadare $[2.].$ A central property of the Joint-Eigenspace Fourier transform and a source of its rich structure is from the proven fact that it factors into the composition of the Poisson transform of the Helgason Fourier transform on $X.$ In this paper we shall employ the Joint-Eigenspace Fourier transform on $X$ to decompose $L^{2}(X)$ into a direct integral and to characterize the irreducibility of the component of the regular representation $T_{X}$ of $G$ on $L^{2}(X).$
	
	\indent Let $X=G/K$ denote a symmetric space of the noncompact type. Let $T_{L^{2}(X)}$ (written as $T_{X}$) and $T_{H_{\lambda}}$ (written as $T_{\lambda}$) denote the regular representations of $G$ on $L^{2}(X)$ and $H_{\lambda},$ respectively; where $H_{\lambda}$ is the image of the Joint-Eigenspace Fourier transform map $f\mapsto(\mathcal{H}_{x}f)(\lambda)$ introduced in Oyadare $[2.]$ as the convolution $(\mathcal{H}_{x}f)(\lambda):=(f\times\varphi_{\lambda})(x),$ $x\in X,$ with $\varphi_{\lambda}$ as the elementary spherical functions on $X.$ A direct integral decomposition of $L^{2}(X)$ reveals the basis nature of the spaces $H_{\lambda}$ and hence of $(\mathcal{H}_{x}f)(\lambda),$ as $f$ runs through $\textbf{D}(X).$ This was used to prove the irreducibility of $T_{\lambda}$ in terms of a pair of members of $\mathfrak{a}^{*}_{\mathbb{C}}.$
	
	\indent The next section gives a summary of basic properties of the Joint-Eigenspace Fourier transform on $X$ (Oyadare $[2.]$) while section three gives the main results.
	
	\ \\
	{\bf \S2. Preliminaries on the Joint-Eigenspace Fourier transform}\\
	({\bf A summary of the results in Oyadare $[2.]$})
	\ \\
	\ \\
	\indent It is well known that a Fourier transform has already been defined and studied on $G/K,$ called the {\it Helgason Fourier transform.} Explicitly, if $f$ is a function defined on $X=G/K,$ then its Helgason Fourier transform $$f\mapsto \hat{f}$$ is given as $$\hat{f}(\lambda,b)=\int_{X}f(x)e^{(-i\lambda+\rho)(A(x,b))}dx$$ for all $\lambda\in \mathfrak{a}^*_{\C},$ $b\in B=K/M$ for which the integral exists $[1.].$ The Helgason Fourier transform is itself an extension (but not a generalization) of the Harish-Chandra Fourier transform as given in the following.
	
	\indent {\bf 2.1 Lemma} (Helgason $[1.],$ p. $224$). The Helgason Fourier transform is an extension of the Harish-Chandra (spherical) Fourier transform from $K-$invariant functions on $X$ to not-necessarily $K-$invariant functions on $X.\;\Box$
	
	\indent This means that the Helgason Fourier transform on $X$ reduces to the formula for the Harish-Chandra spherical Fourier transform when $K-$invariant functions are under consideration.
	
	\indent Just as the spherical-convolution transform on $G$ generalizes the Harish-Chandra Fourier transform on $G$ (which has now been extended to the Helgason Fourier transform on $X$), one would like to seek a generalization of the Helgason Fourier transform on $X$ via the spherical-convolution transform's extension from $G$ to $X=G/K.$ This quest leads us to the consideration of the map $$f\mapsto \mathcal{H}_{x}f$$ for $x\in G/K=X$ defined as $$(\mathcal{H}_{x}f)(\lambda)=(f\times\varphi_{\lambda})(x),$$ where $\varphi_{\lambda}$ is the elementary spherical functions on $X$ and $\times$ is the convolution of functions on $X$ defined as $$(f_{1}\times f_{2})\circ\pi:=(f_{1}\circ\pi)*(f_{2}\circ\pi)$$ for the natural map $\pi:G\rightarrow G/K$ $[1.].$ The following is a cautionary observation in relation to the extensions mentioned above.
	
	\indent {\bf 2.2 Lemma} (Oyadare $[2.],$ p. $10$). The transformation map $f\mapsto \mathcal{H}_{x}f$ defined above for $x\in X,$ $f\in C(X),$ does not restrict to the Helgason Fourier transform on $X.\;\Box$
	
	\indent The above lemma does not erode the fact of Lemma $2.1$ and does not preclude the possibility of the transformation map $f\mapsto \mathcal{H}_{x}f$ coinciding with the Helgason Fourier transform on a special class of functions on $X,$ but it however reveals that we now have a completely different transform on $X$ at hand which we now christen as follows.
	
	\indent {\bf 2.3 Defnition.} Let $x\in X=G/K.$ We shall refer to the transformation map $f\mapsto \mathcal{H}_{x}f$ as the {\it Joint-Eigenspace Fourier transform} on $X$ which is given as $$(\mathcal{H}_{x}f)(\lambda)=(f\times\varphi_{\lambda})(x),$$ for $f\in C_{c}(X),$ $\lambda\in \mathfrak{a}^*_{\C}$ and the elementary spherical functions $\varphi_{\lambda}$ on $X.\;\Box$
	
	\indent For $0<p\leq 2$ and $x\in X,$ the image of the Joint-Eigenspace Fourier transform on $X$ is given as the set $$\mathcal{C}^{p}_{x}(\mathfrak{a}^*_{\C}):=\{\mathcal{H}_{x}f:\lambda\mapsto(\mathcal{H}_{x}f)(\lambda),f\in \mathcal{C}^{p}(X)\}\}.$$ It is thus very important to give explicit description of $\mathcal{C}^{p}_{x}(\mathfrak{a}^*_{\C}).$ The beauty inherent in the Joint-Eigenspace Fourier transform on $X$ (over and above the Helgason Fourier transform) is that the map $f\mapsto\mathcal{H}_{x}f$ could be considered both as a functions of $x$ in $X$ and as a function of $\lambda$ in $\mathfrak{a}^*_{\C}.$
	
	\indent The defining convolution of the Joint-Eigenspace Fourier transform on $X$ was briefly considered  as a function of $X$ by Helgason $[1.]$ (but not as a Fourier transform on $X$) leading to its integral representation in terms of the Helgason Fourier transform on $X$ given as follows.
	
	\indent {\bf 2.4 Lemma} (Helgason $[1.],$ p. $225$). Given that $f\in \mathcal{D}(X):=C^{\infty}_{c}(X)$ and $\lambda\in\mathfrak{a}^*_{\C},$ then $$(\mathcal{H}_{x}f)(\lambda)=\int_{B}e^{(i\lambda+\rho)(A(x,b))}\hat{f}(\lambda,b)db.\;\Box$$
	
	When further considered  as a function of $x,$ Helgason computed the Helgason Fourier transform of the Joint-Eigenspace Fourier transform showing the intimate interconnection between the two Fourier transforms on $X$ which is restated in the notation of Definition $2.3$ as follows.
	
	\indent {\bf 2.5 Lemma} (Oyadare $[2.],$ p. $11$). Let $f\in \mathcal{D}(X).$ Then $$\hat{(\mathcal{H}_{x}f)}(\lambda,b)=\hat{f}(\lambda,b)\cdot\hat{\varphi_{\mu}}(\lambda),$$ $\lambda,\mu\in\mathfrak{a}^*_{\C}.\;\Box$
	
	It is clear that the elementary spherical function $\varphi_{\mu}$ are $K-$bi-invariant as functions on $X.$ In particular, each $\varphi_{\mu}$ is $K-$invariant, so that $\hat{\varphi_{\mu}}(\lambda)=\tilde{\varphi_{\mu}}(\lambda).$ The import of the result of Lemma $2.5$ is that the Helgason Fourier transform of the Joint-Eigenspace Fourier transform (of $f\in C_{c}(X),$ when the Joint-Eigenspace Fourier transform is considered as a function on $X=G/K$) is a non-zero constant multiple (by the constant $c=\tilde{\varphi_{\mu}}(\lambda)$) of the (oridinary) Helgason Fourier transform of $f.$ This connection would help to deduce results for the Joint-Eigenspace Fourier transform from results already established for the Helgason Fourier transform $$f\mapsto\hat{f}(\lambda,b)=\frac{1}{\tilde{\varphi_{\mu}}(\lambda)}\hat{(\mathcal{H}_{x}f)}(\lambda,b).$$
	
	\indent It appears that there was no much interest in any consideration of the Joint-Eigenspace Fourier transform on $X$ as a function of $\lambda\in\mathfrak{a}^*_{\C}$ in much similarity with the detailed treatments of both the Harish-Chandra spherical Fourier transform $f\mapsto\tilde{f}(\lambda)$ and the Helgason Fourier transform $f\mapsto\hat{f}(\lambda,b),$ as such a treatment of the Joint-Eigenspace Fourier transform on $X$ (as a map $f\mapsto(\mathcal{H}_{x}f)(\lambda)$ of $\lambda$) could have led to a better understanding of the joint-eigenspaces on $X.$ Our aim is to kick-start this study. It is however of interest to further note the striking similarity between the inversion formulae for both the Harish-Chandra spherical Fourier transform on $G$ and the Joint-Eigenspace Fourier transform on $X$ established as follows.
	
	\indent {\bf 2.6 Theorem} (Oyadare $[2.],$ p. $12$). For every $f\in \mathcal{D}(X)$ we always have $$f(x)=w^{-1}\int_{\mathfrak{a}^{*}}(\mathcal{H}_{x}f)(\lambda)|c(\lambda)|^{-2}d\lambda,$$ where $w$ is the order of the Weyl group and $c(\lambda)$ is the Harish-Chandra $c-$function$.\;\Box$
	
	\indent The inversion formula for the Joint-Eigenspace Fourier transform on $X,$ given in Theorem $2.6$ above, was curiously noted in $[1.],$ p. $327,$ as the noncompact analogue of the convergent expansion of an arbitrary function on a compact symmetric space and, with the perspective of the Harish-Chandra's inversion formula for the $K-$bi-invariant functions $$F(h)=\int_{K}f(ghkK)dk\;\;\;(f\in C_{c}(X)),$$ we have that $$F(h)=c\int_{\mathfrak{a}^{*}}\hat{F}(\lambda)\varphi_{\lambda}(h)|c(\lambda)|^{-2}d\lambda$$ into which we could insert $\hat{F}(\lambda)=(f\times\varphi_{\lambda})(gK)=(f\times\varphi_{\lambda})(x)=(\mathcal{H}_{x}f)(\lambda).$ Hence Theorem $2.6.$ Our proof in Theorem $2.6$ above is however from the independent recognition of the Joint-Eigenspace Fourier transform as a Fourier transform on $X.$
	
	\indent This now makes it clear that, just as the Helgason Fourier transform is an extension of the Harish-Chandra spherical Fourier transform from $K-$invariant functions on $X$ to not necessarily $K-$invariant functions on $X$ (Recall Lemma $2.1$), the following is the two way bridge between the Helgason Fourier transform on $X$ and the Joint-Eigenspace Fourier transform on $X.$
	
	\indent {\bf 2.7 Theorem} (Oyadare $[2.],$ p. $13$). Let $f\in \mathcal{D}(X)$ and consider the $K-$bi-invariant functions $F(h)=\int_{K}f(ghkK)dk$ in $\textbf{D}(G).$ Then the Joint-Eigenspace Fourier transform on $X$ of any $f\in \mathcal{D}(X)$ is exactly the Helgason Fourier transform of $F$ (on $G$) as an extension of the Harish-Chandra spherical transform to $X.$ That is, $$(\mathcal{H}_{x}f)(\lambda)=\hat{F}(\lambda),\;\;\lambda\in\mathfrak{a}^{*}.\;\Box$$
	
	\indent Theorem $2.6$ above directly leads to the Plancherel formula for the Joint-Eigenspace Fourier transform on $X$ for which we set $$\mathfrak{a}^{*}_{+}=\{\lambda\in\mathfrak{a}:A_{\lambda}\in\mathfrak{a}^{+}\},$$ where $A_{\lambda}$ is characterized by the Killing form requirement that $B(A_{\lambda},H)=\lambda(H)$ for all $H\in\mathfrak{a}.$
	
	\indent {\bf 2.8 Theorem} (Oyadare $[2.],$ p. $14$). The Joint-Eigenspace Fourier transform $f\mapsto(\mathcal{H}_{x}f)(\lambda)$ on $X$ extends to an isometry of $L^{2}(X)$ onto $L^{2}(\mathfrak{a}^{*}_{+}\times B),$ whose measure is given as $|\hat{\varphi_{\mu}}(\lambda)|^{-2}|c(\lambda)|^{-2}d\lambda db$ on $\mathfrak{a}^{*}_{+}\times B.$ Moreover, $$\int_{X}|f(x)|^{2}dx=\int_{\mathfrak{a}^{*}_{+}\times B}|\hat{(\mathcal{H}_{x}f)}(\lambda,b)|^{2}|\hat{\varphi_{\mu}}(\lambda)|^{-2}|c(\lambda)|^{-2}d\lambda db.\;\Box$$
	
	\indent Theorem $2.8$ above gives the group-to-symmetric-space extension of the Plancherel formula for the spherical convolution transform on $G.$ It is observed that we still have on $X$ that $|\hat{\varphi_{\mu}}(\lambda)|^{-2}|c(\lambda)|^{-2}d\lambda db$ is the Plancherel measure of the Joint-Eigenspace Fourier transform on $\mathfrak{a}^{*}_{+}\times B.$ This observation confirms that the $c-$function that is associated with the inversion of the Joint-Eigenspace Fourier transform on $X$ is simply given as $\hat{\varphi_{\mu}}(\lambda)c(\lambda).$
	
	\indent A horocycle in $X=G/K$ is defined to be any orbit $N^{'}\cdot x,$ where $x\in X$ and $N^{'}$ is any subgroup of $G=KAN$ conjugate to $N.$ We shall denote the collection of horocycles on $X$ as $\Xi$ and endow it with the differentiable structure of $G/MN.$ $\Xi$ would be seen as the dual space of $X$) under a very general transform on $X$ defined as follows.
	
	\indent {\bf 2.9 Definition.} The {\it Radon transform} of a function $f$ on $X$ is defined as $\hat{f}(\xi)=\int_{\xi}f(x)ds(x),$ for all $\xi\in\Xi$ for which the integral exists.
	
	\indent The Radon transform is an injective map as explicitly stated below.
	
	\indent {\bf 2.10 Theorem} (Helgason $[1.],$ p. $104$). If $f\in L^{1}(X),$ then $\hat{f}(\xi)$ exists for almost all $\xi\in\Xi$ and if $\hat{f}(\xi)=0$ for almost all $\xi\in\Xi$ then $f(x)=0$ for almost all $x\in X.\;\Box$
	
	\indent Now let $\xi(x,b)$ denote the horocycle passing through the point $x\in X$ with normal $b\in B=K/M.$ We shall denote by $A(x,b)\in\mathfrak{a}$ the composite metric from the origin to $\xi(x,b).$ For $\lambda\in\mathfrak{a}^{*}_{\mathbb{C}},$ $b\in B,$ the function $$e_{\lambda,b}:x\mapsto e_{\lambda,b}(x):=e^{(i\lambda+\rho)(A(x,b))}$$ is a joint-eigenfunction of $\textbf{D}(X)$ and belongs to the joint-eigenspace $\mathcal{E}(X).$ Considering $g\in G$ for which $x=gK\in G/K,$ the integral $$\varphi_{\lambda}(g)=\int_{B}e_{\lambda,b}(gK)=\int_{B}e^{(i\lambda+\rho)(A(gK,b))}$$ is an elementary spherical function on $G$ while $$u_{s}(x)=e_{-is\rho,b}(x)=e^{(s\rho+\rho)(A(x,b))}$$ is a harmonic function on $X$ for each $s\in W.$ We the harmonic function $u_{s=1}(x)=u_{1}(x)=e^{2\rho(A(x,b))}$ is known as the {\it Poisson kernel} ($[6.]$) and this informed the use of $e_{\lambda,b}(x)$ as a kernel to create a map $C(B)\rightarrow\mathcal{E}_{\lambda}(X)$ defined as follows.
	
	\indent {\bf 2.11 Definition.} The {\it Poisson transform} $P_{\lambda}$ of a function $F$ on $B$ is defined as $$P_{\lambda}(x)=\int_{B}e^{(i\lambda+\rho)(A(x,b))}F(b)db,$$ $x\in X.$
	
	\indent The Poisson transform $P_{\lambda}$ maps $C(B)$ into the joint-eigenspace $\mathcal{E}_{\lambda}(X).$ The first connection of the Poisson transform with the Joint-Eigenspace Fourier transform is via Lemma $2.4$ above, giving the Joint-Eigenspace Fourier transform on $X$ as a composition of the Poisson transform on the Helgason Fourier transform on $X$ given as follows.
	
	\indent {\bf 2.12 Theorem} (Oyadare $[2.],$ p. $15$). Given that $f\in \mathcal{D}(X)$ and $\lambda\in\mathfrak{a}^{*}_{\mathbb{C}}$ then $$(\mathcal{H}_{x}f)(\lambda)=(P_{\lambda}\hat{f}(\lambda,\cdot))(x),$$ where $\hat{f}(\lambda,\cdot)\in C(B),$ for each $\lambda.$ In particular, the Joint-Eigenspace Fourier transform on $X$ maps $C(X)$ into the joint-eihenspace $\mathcal{E}_{\lambda}(X).\;\Box$
	
	\indent In other words, the Joint-Eigenspace Fourier transform on $X$ is a Fourier transform on the whole of $X$ (as a lift of the Harish-Chandra spherical Fourier transform on $G$ to the symmetric space $X=G/K$) which does not restrict to the Helgason Fourier transform on $X,$ for all of the entire function space $\mathcal{D}(X).$ It may however be possible to extract a non-empty subspace of $\mathcal{D}(X)$ on which we have the identity $(\mathcal{H}_{x}f)(\lambda)\equiv\hat{f}(\lambda,b),$ for all $x\in X$ and for all $b\in B=K/M.$ Such a subspace of $\mathcal{D}(X),$ if it exists, would contain exactly those $f\in \mathcal{D}(X)$ on which $P_{\lambda}\equiv1,$ the identity map on $C(B).$ Let us therefore set $\mathcal{O}(X)$ to be defined as $$\mathcal{O}(X)=\{f\in\mathcal{D}(X):(\mathcal{H}_{x}f)(\lambda)\equiv\hat{f}(\lambda,b),\;x\in X,\;b\in B=K/M\},$$ whenever it exists.
	
	\indent A summary of the (inter-)relationships among the Helgason Fourier transform, the Poisson transform and the Joint-Eigenspace Fourier transform is presented by the following diagram: $$C(X)\rightarrow^{\hat{f}(\lambda,b)}\rightarrow C(B)\rightarrow^{P_{\lambda}}\rightarrow \mathcal{E}_{\lambda}(X),$$ whose composition is the Joint-Eigenspace Fourier transform on $X.$
	
	\indent The diagram justifies the choice of the name in Definition $2.3$ above. It explains how the Poisson transform completes the duality of the Helgason Fourier transform into a the joint-eigenspace $\mathcal{E}_{\lambda}(X)$ all of which are embodied in the Joint-Eigenspace Fourier transform on $X.$ This observation had also been made for the spherical convolution transform on $G.$ It is therefore clear from Theorem $3.4$ how the surjectivity of the Joint-Eigenspace Fourier transform on $X$ could be deduced from the properties of the Helgason Fourier transform and the Poisson transform.
	
	\indent If, in particular, $\lambda\in\mathfrak{a}^{*}_{\mathbb{C}}$ satisfies $Re(<i\lambda,\alpha>)>0$ for $\alpha\in\sum^{+},$ $b_{o}$ is the origin $eM$ in $B=K/M$ then, for $H\in\mathfrak{a}^{+}$ and $a_{t}=\exp tH,$ we have $$\lim_{t\rightarrow\infty}e^{(-i\lambda+\rho)(tH)}(\mathcal{H}_{a_{t}\cdot o}f)(\lambda)=c(\lambda)\hat{f}(\lambda,b_{o}).$$ Indeed, $$\lim_{t\rightarrow\infty}e^{(-i\lambda+\rho)(tH)}(\mathcal{H}_{a_{t}\cdot o}f)(\lambda)=\lim_{t\rightarrow\infty}e^{(-i\lambda+\rho)(tH)}\int_{B}e^{(i\lambda+\rho)(A(a_{t}\cdot o,b))}\hat{f}(\lambda,b)db$$ $$=c(\lambda)\hat{f}(\lambda,b_{o}),$$ by Theorem $3.16$ of $[1.],$ p. $120.$
	
	\indent {\bf 2.13 Lemma} (Oyadare $[2.],$ p. $17$). Given that $Re(<i\lambda,\alpha>)>0$ for $\alpha\in\sum^{+}$ and some $\lambda\in\mathfrak{a}^{*}_{\mathbb{C}},$ then, for $a_{t}=\exp tH,$ $H\in\mathfrak{a}^{+}$ and $b_{o}$ as the origin $eM$ in $B=K/M,$ we have that $$\lim_{t\rightarrow\infty}e^{(-i\lambda+\rho)(tH)}(\mathcal{H}_{a_{t}\cdot o}f)(\lambda)=c(\lambda)\hat{f}(\lambda,b_{o}).\;\Box$$
	
	\indent We equally have a functional equation for the Joint-Eigenspace Fourier transform on $X$ given as $$\int_{K}(\mathcal{H}_{gk\cdot x}f)(\lambda)dk=\varphi_{\lambda}(x)(\mathcal{H}_{g\cdot \bar{e}}f)(\lambda)$$ for $x\in X,$ $g\in G,$ $f\in C^{\infty}_{c}(X).$
	
	\indent In order to now establish a Paley-Wiener theorem for the Joint-Eigenspace Fourier transform on $X$ we shall examine closely both the Paley-Wiener theorem for the Helgason Fourier transform on $X$ and the bijectivity of the Poisson transform. To this end, we shall refer to a $C^{\infty}$ function $\psi(\lambda,b),$ $\lambda\in\mathfrak{a}^{*}_{\mathbb{C}},$ $b\in B,$ as being a {\it holomorphic function of uniform exponential type} and hence belonging to the space $\mathcal{H}^{R}(\mathfrak{a}^{*}\times B),$ if it is holomorphic in $\lambda$ and if there exists a constant $R\geq0$ such that for each $N\in\mathbb{Z}^{+},$ $$\sup_{\lambda\in\mathfrak{a}^{*}_{\mathbb{C}},b\in B}e^{-R|Im \lambda|}(1+|\lambda|)^{N}|\psi(\lambda,b)|<\infty,$$ where $Im \lambda$ is the imaginary part of $\lambda$ and $|\lambda|$ is its modulus. We set $$\mathcal{H}(\mathfrak{a}^{*}\times B)=\bigcup_{R>0}\mathcal{H}^{R}(\mathfrak{a}^{*}\times B)$$ and denote by $\mathcal{H}(\mathfrak{a}^{*}\times B)_{W}$ as members of $\mathcal{H}(\mathfrak{a}^{*}\times B)$ which are Weyl group invariant. Members of $\mathcal{H}(\mathfrak{a}^{*}\times B)_{W}$ may be viewed as being of the form $\psi(\lambda,b)=\psi_{\lambda}(b),$ so that $\mathcal{H}(\mathfrak{a}^{*}\times B)_{W}\subset C(B)=$ the domain of the Poisson transform, $P_{\lambda}.$ Thus the Weyl group invariance of members $\psi(\lambda,b)=\psi_{\lambda}(b)\in\mathcal{H}(\mathfrak{a}^{*}\times B)_{W}$ may be stated in terms of $P_{\lambda}$ as $P_{s\lambda}(\psi_{s\lambda})=P_{\lambda}(\psi_{\lambda})$ for each $s\in W.$ That is, $$\int_{B}e^{(is\lambda+\rho)(A(x,b))}\psi_{s\lambda}(b)db=\int_{B}e^{(i\lambda+\rho)(A(x,b))}\psi_{\lambda}(b)db$$ for each $s\in W.$ The following is the crucial Paley-Wiener theorem for the Helgason Fourier transform.
	
	\indent {\bf 2.14 Theorem} (Helgason $[1.],$ p. $270$). The Helgason Fourier transform $f\mapsto\hat{f}(\lambda,b)$ is a bijection of $C^{\infty}_{c}(X)$ onto $\mathcal{H}(\mathfrak{a}^{*}\times B)_{W}.$ Moreover, $\psi=\hat{f}$ is in $\mathcal{H}(\mathfrak{a}^{*}\times B)_{W}$ iff $supp(f)\subset Cl(B_{R}(0)).\;\Box$
	
	\indent The earlier diagram connecting the three transforms of this paper now becomes refined as $$C^{\infty}_{c}(X)\rightarrow^{\hat{f}(\lambda,b)}\rightarrow \mathcal{H}(\mathfrak{a}^{*}\times B)_{W}\rightarrow^{P_{\lambda}}\rightarrow \mathcal{E}_{\lambda}(X).$$
	
	\indent It is now very clear that it is not a coincidence that the Poisson transform is naturally included in the definition of the Paley-Wiener space for the Helgason Fourier transform on $X.$ Our contribution in this respect is to explore this fact to complement the journey of the Helgason Fourier transform on $X.$ We believe that this natural involvement of the Poisson transform in the image-construction of the Paley-Wiener space for the Helgason Fourier transform on $X$ is meant to complement the Helgason Fourier transform and to build it up into the status of a Joint-Eigenspace Fourier transform on $X.$ We now embark on the completion of the building up of Helgason Fourier transform on $X.$
	
	\indent Let $e^{\beta}$ be the $e_{s^{*}}-$function for $G_{\beta}/K_{\beta}$ in which $$e_{s^{*}}(\lambda)=const\cdot\prod_{\beta\in\sum^{+}_{o}}e^{\beta}(\lambda_{\beta}).$$ We know, from $[1.],$ p. $269,$ that $P_{\lambda}$ is injective iff $e_{s^{*}}(\lambda)\neq0.$ On the surjectivity of $P_{\lambda}$ we note that if $\nu:\textbf{D}(G)\rightarrow\textbf{E}(X)$ denote the homomorphism induced by the action $f\mapsto f^{\tau(g)}$ of $G$ on $\mathcal{E}(X)$ and if $$\mathcal{E}^{\infty}_{\lambda}(X)=\{f\in\mathcal{E}_{\lambda}(X): (\nu(D)f)(x)=O(e^{A\;d(0,x)}),\forall D\in\textbf{D}(G)\;\mbox{and some}\;A>0\},$$ then for each $\lambda\in\mathfrak{a}^{*}_{\mathbb{C}}$ for which $e_{s^{*}}(\lambda)\neq0,$ we have that $$P_{\lambda}(\mathcal{E}(B))=\mathcal{E}^{\infty}_{\lambda}(X).$$
	
	\indent Since the $\mathcal{H}(\mathfrak{a}^{*}\times B)_{W}\subset C(B)$ we can now restrict $P_{\lambda}$ to $\mathcal{H}(\mathfrak{a}^{*}\times B)_{W}.$ Clearly $P_{\lambda}(\mathcal{H}(\mathfrak{a}^{*}\times B)_{W})\subset\mathcal{E}^{\infty}_{\lambda}(X)$ and if we define the Hilbert space $H^{\infty}_{\lambda}(X)$ as $$H^{\infty}_{\lambda}(X):=H_{\lambda}(X)\bigcap\mathcal{E}^{\infty}_{\lambda}(X),$$ from the Hilbert space $$H_{\lambda}(X):=\{h\in\mathcal{E}_{\lambda}(X):h(x)=\int_{B}e^{(i\lambda+\rho)(A(x,b))}F(b),\;F\in L^{2}(B)\}$$ ($[1.],$ p. $531$), then $$P_{\lambda}(\mathcal{H}(\mathfrak{a}^{*}\times B)_{W})=H^{\infty}_{\lambda}(X),$$ where $\lambda\in\mathfrak{a}^{*}_{\mathbb{C}}$ is simple. Thus a Paley-Wiener theorem for the Joint-Eigenspace Fourier transform is immediate and is as given in the following main result of the paper.
	
	\indent {\bf 2.15 Theorem} (Oyadare $[2.],$ p. $19$). Let $\lambda\in\mathfrak{a}^{*}_{\mathbb{C}}$ be simple. The Joint-Eigenspace Fourier transform on $X$ is a bijection of $C^{\infty}_{c}(X)$ onto $H^{\infty}_{\lambda}(X).$ Moreover, we have that $\psi(x)=(\mathcal{H}_{x}f)(\lambda)$ is in $H^{\infty}_{\lambda}(X)$ iff $supp(f)\subset Cl(B_{R}(0)).$ $\;\Box$
	
	\indent We are now in a position to characterize the earlier defined $\mathcal{O}(X)-$subspace of the domain $\mathcal{D}(X)$ of the Joint-Eigenspace Fourier transform on $X.$
	
	\indent {\bf 2.16 Corollary} (Oyadare $[2.],$ p. $20$). Let $\lambda\in\mathfrak{a}^{*}_{\mathbb{C}}$ be simple. Then the subspace $\mathcal{O}(X)$ exists iff $\mathcal{H}(\mathfrak{a}^{*}\times B)_{W}=H^{\infty}_{\lambda}(X).\;\Box$
	
	\indent The Paley-Wiener version of our diagram is now commutative and is finally given (for $\lambda\in\mathfrak{a}^{*}_{\mathbb{C}}$ simple) as 
	
	$$C^{\infty}_{c}(X)\rightarrow^{\hat{f}(\lambda,b)}\rightarrow \mathcal{H}(\mathfrak{a}^{*}\times B)_{W}\rightarrow^{P_{\lambda}}\rightarrow \mathcal{E}_{\lambda}(X).$$
	
	\ \\
	{\bf \S3. $L^{2}-$Theory}
	\ \\
	\ \\
	\indent We recall the image $H_{\lambda}(X)$ of the Joint-Eigenspace Fourier transform $f\mapsto(\mathcal{H}_{x}f)(\lambda)$ for $C^{\infty}_{c}(X),$ given as $$H_{\lambda}(X):=\{h\in\mathcal{E}_{\lambda}(X):h(x)=\int_{B}e^{(i\lambda+\rho)(A(x,b))}F(b)db,\;F\in L^{2}(B)\}.$$ It follows therefore that the Fourier transform of $f\in C^{\infty}_{c}(X)$ in $H_{\lambda}(X)$ is $(\mathcal{H}_{x}f)(\lambda)$ as a function of $x.$ We already know that for $\lambda\in\mathfrak{a}^{*}_{\mathbb{C}}$ simple, $H_{\lambda}(X)$ becomes a Hilbert space under the $\lambda-$norm, $|\cdot|_{\lambda},$ defined as $|h|_{\lambda}:=||F||_{L^{2}(B)}.$ This means that we now have the regular representations $T_{L^{2}(X)}$ $(=:T_{X})$ and $T_{\mathcal{E}_{\lambda}(X)}$ $(=:T_{\lambda})$ of $G$ on the Hilbert spaces $L^{2}(X)$ and $\mathcal{E}_{\lambda}(X)$ given respectively as $$(T_{X}(g)f)(x):=f(g^{-1}\cdot x),\;\;x\in X,g\in G, f\in L^{2}(X)$$ and $$(T_{\lambda}(g)f)(x):=f(g^{-1}\cdot x),\;\;x\in X,g\in G, f\in \mathcal{E}_{\lambda}(X).$$ The Joint-Eigenspace transform then gives a direct integral decomposition of $L^{2}(X)$ and $T_{X}$ as follows.
	
	\indent {\bf 3.1 Theorem} With the notation as above, the Joint-Eigenspace Fourier transform gives $$L^{2}(X)=\int_{\mathfrak{a}^{*}/W}H_{\lambda}(X)|c(\lambda)|^{-2}d\lambda$$ and $$T_{X}=\int_{\mathfrak{a}^{*}/W}T_{\lambda}|c(\lambda)|^{-2}d\lambda$$
	
	\indent {\bf Proof.} For $f_{1},f_{2}\in L^{2}(X),$ $$\int_{X}f_{1}(x)\overline{f_{2}(x)}dx=\frac{1}{\omega}\int_{\mathfrak{a}^{*}\times B}\hat{f_{1}}(\lambda,b)\overline{\hat{f_{2}}(\lambda,b)}|c(\lambda)|^{-2}d\lambda db$$ $$=\frac{1}{\omega}\int_{\mathfrak{a}^{*}\times B}[\hat{(\mathcal{H}_{x}f_{1})}(\lambda,b)\hat{\varphi_{\mu}}(\lambda)^{-1}]
	\overline{[\hat{(\mathcal{H}_{x}f_{2})}(\lambda,b)\hat{\varphi_{\mu}}(\lambda)^{-1}]}|c(\lambda)|^{-2}d\lambda db$$
	$$=\frac{1}{\omega}\int_{\mathfrak{a}^{*}\times B}\hat{(\mathcal{H}_{x}f_{1})}(\lambda,b)\overline{\hat{(\mathcal{H}_{x}f_{2})}(\lambda,b)}
	\hat{\varphi_{\mu}}(\lambda)^{-1}
	\overline{\hat{\varphi_{\mu}}(\lambda)^{-1}}|c(\lambda)|^{-2}d\lambda db$$
	$$=\frac{1}{\omega}\int_{\mathfrak{a}^{*}\times B}\hat{(\mathcal{H}_{x}f_{1})}(\lambda,b)\overline{\hat{(\mathcal{H}_{x}f_{2})}(\lambda,b)}
	|\hat{\varphi_{\mu}}(\lambda)|^{-2}|c(\lambda)|^{-2}d\lambda db.$$ With $f_{1}=f_{2},$ we then have that $$\int_{X}|f(x)|^{2}dx=\int_{X}f(x)\overline{f(x)}dx$$ $$=\frac{1}{\omega}\int_{\mathfrak{a}^{*}\times B}\hat{(\mathcal{H}_{x}f)}(\lambda,b)\overline{\hat{(\mathcal{H}_{x}f)}(\lambda,b)}
	|\hat{\varphi_{\mu}}(\lambda)|^{-2}|c(\lambda)|^{-2}d\lambda db$$
	$$=\frac{1}{\omega}\int_{\mathfrak{a}^{*}\times B}|\hat{(\mathcal{H}_{x}f_{1})}(\lambda,b)|^{2}
	|\hat{\varphi_{\mu}}(\lambda)|^{-2}|c(\lambda)|^{-2}d\lambda db,$$ in which the Plancherel measure $\omega^{-1}|\hat{\varphi_{\mu}}(\lambda)|^{-2}|c(\lambda)|^{-2}d\lambda db$ on $\mathfrak{a}^{*}_{+}\times B$ is seen as a non-zero multiple of the well-known Harish-Chandra Plancherel measure $\omega^{-1}|c(\lambda)|^{-2}d\lambda db.\;\Box$
	
	\indent The last Theorem reveals the basis nature of the $H_{\lambda}-$spaces with respect to the Hilbert space $L^{2}(X)$ and hence the basis nature of the Joint-Eigenspace Fourier transform $f\mapsto(\mathcal{H}_{x}f)(\lambda)$ with respect to both Hilbert spaces $H_{\lambda}$ and $L^{2}(X).$ These decompositions are clearly not available for the Helgason Fourier transform on $X.$ The reader may also consult Helgason $[1.],$ p. $546,$ for some results which {\it have now been put in their proper perspectives with respect to the Joint-Eigenspace Fourier transform.}
	
	\indent This theorem also verifies the reducibility of $T_{X}$ on $L^{2}(X)$ via the Joint-Eigenspace Fourier transform. In order to make the above direct integral decompositions very useful it is now crucial to study the irreducbility of $T_{\lambda}$ on $H_{\lambda}(X).$
	
	\indent If $-\lambda\in\mathfrak{a}^{*}_{\mathbb{C}}$ is simple, then the space of functions $b\mapsto \hat{f}(\lambda,b)$ is dense in $L^{2}(B)$ as $f$ runs through $C^{\infty}_{c}(X),$ $[1.],$ p. $227.$ In other words, the space $$C^{\infty,\lambda}_{c}(X):=\{\hat{f_{\lambda}}:b\mapsto \hat{f_{\lambda}}(b):=\hat{f}(\lambda,b),\;f\in C^{\infty}_{c}(X)\}$$ is dense in $L^{2}(B),$ for every simple $-\lambda\in\mathfrak{a}^{*}_{\mathbb{C}}.$ This translates to mean that every $F\in L^{2}(B)$ is of the form $$F(b)=\sum_{k}a_{k}\hat{f}_{k}(\lambda,b),$$ for simple $-\lambda\in\mathfrak{a}^{*}_{\mathbb{C}},$ $a_{k}\in\mathbb{C}.$
	
	\indent {\bf 3.2 Lemma} Let $-\lambda\in\mathfrak{a}^{*}_{\mathbb{C}}$ be simple. The closed subspace $\mathcal{E}_{(\lambda)}(X)$ of $\mathcal{E}(X)$ generated by the $G-$translates of $\varphi_{\lambda}$ contains all of $(\mathcal{H}_{x}f)(\lambda)$ and its linear combinations, as $f$ runs through $C^{\infty}_{c}(X).$
	
	\indent {\bf Proof.} The subspace $\mathcal{E}_{(\lambda)}(X)$ contains all of $f$ given as $$f(x)=\int_{B}e^{(i\lambda+\rho)(A(x,b))}F(b)db$$ where $F\in L^{2}(B),$ $[1.],$ p. $234.$ Since every $F\in L^{2}(B)$ is of the form $$F(b)=\sum_{k}a_{k}\hat{f}_{k}(\lambda,b),$$ with $-\lambda\in\mathfrak{a}^{*}_{\mathbb{C}},$ simple, $a_{k}\in\mathbb{C},$ then $$f(x)=\int_{B}e^{(i\lambda+\rho)(A(x,b))}F(b)db$$ $$=\int_{B}e^{(i\lambda+\rho)(A(x,b))}(\sum_{k}a_{k}\hat{f}_{k}(\lambda,b))db=\sum_{k}a_{k}(\int_{B}e^{(i\lambda+\rho)(A(x,b))}\hat{f}_{k}(\lambda,b)db)$$
	$$=\sum_{k}a_{k}((\mathcal{H}_{x}f_{k})(\lambda)).\;\Box$$
	
	\indent It is known that the irreducibility question for the regular representation $T_{\lambda}$ of $G$ on $H_{\lambda}$ amounts to the problem of finding for which $\lambda\in\mathfrak{a}^{*}_{\mathbb{C}}$ the spaces $\mathcal{E}_{(\lambda)}(X)$ and $\mathcal{E}_{\lambda}(X)$ coincide, $ [1.],$ p. $526.$ This problem is addressed in the following results.
	
	\indent {\bf 3.3 Theorem.} $\mathcal{E}_{(\lambda)}(X)=\mathcal{E}_{\lambda}(X)$ iff $\lambda$ and $-\lambda$ are both simple in $\mathfrak{a}^{*}_{\mathbb{C}}.$
	
	\indent {\bf Proof.} Let both $\lambda$ and $-\lambda$ be simple in $\mathfrak{a}^{*}_{\mathbb{C}}.$ Then $\mathcal{E}_{(\lambda)}(X)\subseteq\mathcal{E}_{\lambda}(X),$ since $D\varphi_{\lambda}=\Gamma(D)(i\lambda)\varphi_{\lambda},$ for every $D\in\textbf{D}(X).$ Also $\mathcal{E}_{\lambda}(X)\subseteq\mathcal{E}_{(\lambda)}(X),$ since $L^{2}(X)=\int_{\mathfrak{a}^{*}/W}H_{\lambda}(X)|c(\lambda)|^{-2}d\lambda,$ in which $$H_{\lambda}=\mbox{image of the Joint-Eigenspace Fourier transform}\;f\mapsto(\mathcal{H}_{x}f)(\lambda)$$ is contained in $\mathcal{E}_{(\lambda)}(X),$ as established in Lemma $3.2.$
	
	\indent Conversely, if $\mathcal{E}_{(\lambda)}(X)=\mathcal{E}_{\lambda}(X)$ for some $\lambda\in\mathfrak{a}^{*}_{\mathbb{C}},$ where $H_{\lambda}$ is the aforementioned dense subspace of both, then the Joint-Eigenspace Fourier transform $f\mapsto(\mathcal{H}_{x}f)(\lambda)$ is a well-defined map into $H_{\lambda}.$ Now as $(\mathcal{H}_{x}f)(\lambda)=(P_{\lambda}\hat{f}(\lambda,\cdot))(x),$ the injectivity of the Poisson transform $P_{\lambda}$ implies that $\lambda$ is simple in $\mathfrak{a}^{*}_{\mathbb{C}}.$ Also $\mathcal{E}_{(\lambda)}(X)=\mathcal{E}_{\lambda}(X)$ implies that $\mathcal{E}_{(-\lambda)}(X)=\mathcal{E}_{-\lambda}(X),$ which in turn implies that $-\lambda\in\mathfrak{a}^{*}_{\mathbb{C}}$ is simple$.\;\Box$
	
	\indent It is known that each $H_{\lambda}$ is dense in $\mathcal{E}_{\lambda}(X)$ (for simple $\lambda\in\mathfrak{a}^{*}_{\mathbb{C}}$) $[1.],$ p. $296.$ The fact that $\mathcal{E}_{(\lambda)}(X)$ contains all of $(\mathcal{H}_{x}f)(\lambda)$ (for simple $-\lambda\in\mathfrak{a}^{*}_{\mathbb{C}}$) implies that $\mathcal{E}_{(\lambda)}(X)$ also contains the dense subspace $H_{\lambda}.$ Thus we have the following.
	
	\indent {\bf 3.4 Lemma.} $\mathcal{E}_{(\lambda)}(X)=span(H_{\lambda})=\mathcal{E}_{\lambda}(X)$ for $\lambda$ and $-\lambda$ both simple in $\mathfrak{a}^{*}_{\mathbb{C}}.\;\Box$
	
	\indent We are now in a position to state our main result as follows.
	
	\indent {\bf 3.5 Theorem.} $T_{\lambda}$ is irreducible iff $\lambda$ and $-\lambda$ are both simple in $\mathfrak{a}^{*}_{\mathbb{C}}.$
	
	\indent {\bf Proof.} $T_{\lambda}$ is irreducible iff $\mathcal{E}_{(\lambda)}(X)=\mathcal{E}_{\lambda}(X)$ iff $\lambda$ and $-\lambda$ are both simple in $\mathfrak{a}^{*}_{\mathbb{C}}.\;\Box$
	
	\ \\
	{\bf \S4. Conclusion.}
	
	This paper is a short note on characterization of the irreducibility of the regular representation $T_{\lambda}$ of $G$ on $H_{\lambda}.$ It further proves the rich structure of the Joint-Eigenspace Fourier transform on $G/K$ and its advantage over the Helgason Fourier transform. The present paper is a testament to the synergy brought into the geometric analysis of noncompact symmetric spaces by the Joint-Eigenspace Fourier transform of Oyadare $[2.].$ The results of the present paper may be used for the classification of some distinguished irreducible symmetric spaces $[3.].$
	\ \\
	\ \\
	{\bf Declarations.}
	\ \\
	
	{\it Ethics Approval and consent to participate:} I declare that there is no unethical conduct in the research and I also give my consent to participate in Ethics Approval.
	\ \\
	
	{\it Consent for publication:} I give my consent to the publication of this manuscript.
	\ \\
	
	{\it Availability of data and materials:} Materials consulted during the research are as listed in the References below.
	\ \\
	
	{\it Competing interest:} I declare that there is no competing interest in the conduct of this research.
	\ \\
	
	{\it Funding:} I declare that I received no funding for this research.
	\ \\
	
	{\it Authors' contributions:} I declare that I was the only one who contributed to the conduct of this research.
	\ \\
	
	{\it Acknowledgments:} I acknowledge authors whose research works were consulted as listed in the References below.
	\ \\
	\ \\
	\ \\
	{\bf   References.}
	\begin{description}
		
\item [{[1.]}]  Helgason, S., \textit{Geometric analysis on symmetric spaces,} Mathematical Surveys and Monographs, vol. $39,$ Providence, Rhode Island $(1994)$

\item [{[2.]}] Oyadare, O. O.,  A Paley-Wiener theorem for the Joint-Eigenspace Fourier transform on noncompact symmetric spaces, {\it arXiv:$2409.09036.$} [math.FA] $(2024)$

\item [{[3.]}] Tojo, K., Classification of irreducible symmetric spaces which admit standard compact Clifford-Klein form, {\it Proc. Japan Acad.,} {\bf 95}. Ser. $A$ $(2019)$
	\end{description}
\end{document}